\documentclass[10pt,a4paper,leqno]{amsart}
\makeatletter
\renewcommand\subsection{\@startsection{subsection}{2}%
	\z@{.5\linespacing\@plus.7\linespacing}{-.5em}%
	{\normalfont\scshape}}
\makeatother
\usepackage{hyperref}
\usepackage{amssymb,amsmath,amsthm,fancyhdr,graphicx,MnSymbol,latexsym}

\newcommand{\eqand}{\ensuremath{\quad \textrm{and} \quad}}

\newcommand{\foot}{\footnote}

\newcommand{\kl}{\left(}

\newcommand{\ml}{\left\{}
\newcommand{\il}{\left[}

\newcommand{\mm}{\,\middle|\,}

\newcommand{\kr}{\right)}

\newcommand{\mr}{\right\}}
\newcommand{\ir}{\right]}

\newcommand{\ld}{\ensuremath{,\ldots,}}
\newcommand{\ssq}{\ensuremath{\subseteq}}
\newcommand{\smin}{\ensuremath{\setminus}}
\newcommand{\eps}{\ensuremath{\varepsilon}}

\newcommand{\ind}{\ensuremath{\mathbf{1}}}

\newcommand{\supp}{\ensuremath{\mathrm{supp}}}

\newcommand{\diam}{\ensuremath{\mathrm{diam}}}

\newcommand{\alphlist}{\begin{list}{(\alph{enumi})}{\usecounter{enumi}\setlength{\parsep}{2pt}
      \setlength{\itemsep}{1pt} \setlength{\topsep}{5pt}
      \setlength{\partopsep}{3pt}}}
\newcommand{\arablist}{\begin{list}{(\arabic{enumi})}{\usecounter{enumi}\setlength{\parsep}{2pt}
          \setlength{\itemsep}{1pt} \setlength{\topsep}{5pt}
          \setlength{\partopsep}{3pt}}}
\newcommand{\romanlist}{\begin{list}{(\roman{enumi})}{\usecounter{enumi}\setlength{\parsep}{2pt}
              \setlength{\itemsep}{1pt} \setlength{\topsep}{5pt}
              \setlength{\partopsep}{3pt}}}
\newcommand{\Romanlist}{\begin{list}{(\Roman{enumi})}{\usecounter{enumi}\setlength{\parsep}{2pt}
              \setlength{\itemsep}{1pt} \setlength{\topsep}{5pt}
              \setlength{\partopsep}{3pt}}}
\newcommand{\bulletlist}{\begin{list}{$\bullet$}{\setlength{\parsep}{2pt}
                \setlength{\itemsep}{1pt} \setlength{\topsep}{5pt}
                \setlength{\partopsep}{3pt}\setlength{\leftmargin}{15pt}}} 
\newcommand{\Alphlist}{\begin{list}{(\Alph{enumi})}{\usecounter{enumi}\setlength{\parsep}{2pt}
      \setlength{\itemsep}{1pt} \setlength{\topsep}{5pt}
      \setlength{\partopsep}{3pt}}}
 \newcommand{\listend}{\end{list}}

\newcommand{\N}{\ensuremath{\mathbb{N}}} 
\newcommand{\R}{\ensuremath{\mathbb{R}}}
\newcommand{\Z}{\ensuremath{\mathbb{Z}}}

\newcommand{\cA}{\mathcal{A}}
\newcommand{\cB}{\mathcal{B}}
\newcommand{\cC}{\mathcal{C}}

\newcommand{\cU}{\mathcal{U}}

\newcommand{\nLim}{\ensuremath{\lim_{n\rightarrow\infty}}}

\newcommand{\inergsum}{\ensuremath{\sum_{i=0}^{n-1}}}

\newcommand{\ktel}{\ensuremath{\frac{1}{k}}}

\newcommand{\ntel}{\ensuremath{\frac{1}{n}}}

\usepackage{papercommands}

\newcommand{\Dmmax}{\ensuremath{\mathrm{D}_m^{\max}}}
\newcommand{\Dm}{\ensuremath{\mathrm{D}_m}}
\newcommand{\Dmbar}{\ensuremath{\overline{\mathrm{D}}_m^\varphi}}
\newcommand{\Dmymax}{\ensuremath{\mathrm{D}_{Y,m}^{\max}}}

\newcommand{\Dmimax}{\ensuremath{\mathrm{D}_{m+1}^{\max}}}
\newcommand{\Dmi}{\ensuremath{\mathrm{D}_{m+1}}}
\newcommand{\Dmibar}{\ensuremath{\overline{\mathrm{D}}_{m+1}^\varphi}}
\newcommand{\Dmiymax}{\ensuremath{\mathrm{D}_{Y,m+1}^{\max}}}

\title[Multivariate mean equicontinuity for finite-to-one topomorphic
  extensions] {Multivariate mean equicontinuity for finite-to-one topomorphic
  extensions} \author{J.~Breitenb\"ucher \and L.~Haupt \and
  T.~J\"ager}

\begin{document}

         \begin{abstract}In this note, we generalise the concept of topo-isomorphic
          extensions and define finite topomorphic extensions as topological
          dynamical systems whose factor map to the maximal equicontinuous factor
          is measure-theoretically at most $m$-to-one for some $m\in\N$. We
          further define multivariate versions of mean equicontinuity,
          complementing the notion of multivariate mean sensitivity introduced
          by Li, Ye and Yu, and then show that any $m$-to-one topomorphic
          extension is mean $(m+1)$-equicontinuous. This falls in line with the
          well-known result, due to Downarowicz and Glasner, that strictly
          ergodic systems are isomorphic extensions if and only if they are mean
          equicontinuous. While in the multivariate case we can only conjecture
          that the converse direction also holds, the result provides an
          indication that multivariate equicontinuity properties are strongly
          related to finite extension structures. For minimal systems, an
          Auslander-Yorke type dichotomy between multivariate mean
          equicontinuity and sensitivity is shown as well.
         \end{abstract}

\maketitle


\section{Introduction}

A topological system (tds) $(X,\varphi)$, given by a compact metric space $X$
and a homeomorphism $\varphi:X\to X$, is called {\em mean equicontinuous} if the
{\em Besicovitch pseudo-metric} $d_B$ given by
\[
d_B(x,y) \ = \ \limsup_{n\to\infty} \ntel \inergsum
d\kl\varphi^i(x),\varphi^i(y)\kr
\]
is continuous with respect to the original metric $d$ on $X$. This notion was
introduced by Li, Tu and Ye in \cite{LiTuYe2015MeanSensitivity}, who also proved
that any minimal mean equicontinuous tds is uniquely ergodic. Moreover,
Downarowicz and Glasner showed that this property is closely related to the
extension structure of $(X,\varphi)$. We denote by $(Y,\psi)$ the {\em maximal
  equicontinuous factor} (MEF) of $(X,\varphi)$ and by $\pi: X\to Y$ the
corresponding factor map. Then in the minimal case $(X,\varphi)$ is mean
equicontinuous if and only if it is uniquely ergodic and $\pi$ is a
measure-theoretic isomorphism between the two systems $(X,\varphi)$ and $(Y,\psi)$
equipped with their respective unique invariant probability measures
\cite{DownarowiczGlasner2015IsomorphicExtensionsAndMeanEquicontinuity}. In this
situation, we say $(X,\varphi)$ is a {\em topo-isomorphic extension} of
$(Y,\psi)$.  This seminal result prompted further research in different
directions. It was generalised in
\cite{FuhrmannGroegerLenz2022MeanEquicontinuousGroupActions} to the non-minimal
case and more general group actions. In
\cite{GarciaJaegerYe2021DiamMeanEquicontinuity}, various subclasses of
topo-isomorphic extensions, defined in terms of additional invertibility
properties of the factor map $\pi$, were characterised by intrinsic dynamical
properties of the system. At the same time, multivariate versions of mean
sensitivity -- the counterpart to mean equicontinuity -- were introduced by Li,
Ye and Yu \cite{LiYeYu2022EquicontinuityAndSensitivityInMeanForms} and further
studied by various authors
\cite{LiYu2021MeanSensitiveTuples,LiuYin2023MeanSensitiveTuplesOfGroupActions}.

Broadly speaking, the aim of this note is to establish a link between the
multivariate version of mean equicontinuity -- complementary to the notion of
multivariate mean sensitivity in
\cite{LiYeYu2022EquicontinuityAndSensitivityInMeanForms} -- and the extension
structure of the system, similar to the result on topo-isomorphic extensions in
\cite{DownarowiczGlasner2015IsomorphicExtensionsAndMeanEquicontinuity}. We also
refer to
\cite{ShaoYeZhang2008NSensitivity,HuangLianShaoYe2021MinimalSystemWithFinitelyManyMeasures}
for related results concerning the interplay between dynamical properties
finite-to-one extension structures.

In order to be more precise, we need to introduce some notation. Given $m\in\N,
m\geq 2$ and $x_1\ld x_m\in X$, we let
\[
\Dm(x_1\ld x_m)\ =\ \min_{1\leq i<j\leq m} d(x_i,x_j) \eqand
\Dmmax(x_1\ld x_m)\ = \ \max_{1\leq i<j\leq m} d(x_i,x_j) \ .
\]
We define the {\em Besicovitch $m$-distance} of $x_1\ld x_m$ as
\[
	\Dmbar(x_1\ld x_m)\ = \ \limsup_{n\to\infty} \ntel \inergsum
        \Dm\kl\varphi^i(x_1)\ld \varphi^i(x_m)\kr \
\]
and say $(X,\varphi)$ is {\em mean $m$-equicontinuous} if for all $\eps>0$ there
exists $\delta>0$ such that $\Dmmax(x_1\ld x_m)<\delta$ implies $\Dmbar(x_1\ld
x_m)<\eps$.

The corresponding extension structure is defined as follows. First, given a map
$\gamma:Y\to X^m$ that satisfies $\pi\circ\gamma_i(y)=y$ for all $y\in Y$ and $i \in \ml 1, \ldots, n\mr$, we
denote the corresponding point set by $\Gamma=\{\gamma_i(y)\mid y\in Y,\ i=1\ld
m\}$. Then, we call $(X,\varphi)$ an {\em $m$:1 topomorphic extension (of its
  MEF)} if $m$ is the least integer such that there exists a measurable map
$\gamma:Y\to X^n$ which satisfies $\mu(\Gamma)=1$ for any $\varphi$-invariant
measure $\mu$ on $X$. Equivalently, one could require that any
$\varphi$-invariant measure is supported on at most $m$ points in every fibre
$\pi^{-1}(y),\ y\in Y$ (see Section~\ref{FTE} for details and further
discussion). We note that mean $2$-equicontinuity is just mean equicontinuity
and a 1:1 topomorphic extension is just a topo-isomorphic extension in the sense
of
\cite{LiTuYe2015MeanSensitivity,DownarowiczGlasner2015IsomorphicExtensionsAndMeanEquicontinuity,FuhrmannGroegerLenz2022MeanEquicontinuousGroupActions}. Hence,
for minimal systems mean 2-equicontinuity and a 1:1 topomorphic extension
structure are equivalent by
\cite{DownarowiczGlasner2015IsomorphicExtensionsAndMeanEquicontinuity}, as
mentioned above.\medskip

Here, our main result reads as follows.
\begin{thm}\label{t.multivariate_mean_equicont_for_topomorphic_extensions}
  Let $m\in\N$ and suppose that $(X,\varphi)$ is a minimal $m$:1 topomorphic
  extension of its MEF $(Y,\psi)$. Then $(X,\varphi)$ is mean
  $(m+1)$-equicontinuous.
\end{thm}
We actually believe that, as in the case $m=1$, the converse holds as
well. However, we do not pursue this problem here and only discuss briefly in
Section~\ref{MainResult} why a proof of this fact -- if it is true -- needs to
be more intricate in the multivariate case.
\begin{conj} \label{conj.topomorphic_extensions}
   Let $m\in\N$ and suppose that $(X,\varphi)$ is minimal. Then $(X,\varphi)$ is
   mean \mbox{$(m+1)$-equicontinuous,} but not mean $m$-equicontinuous, if and
   only if it is an $m$:1 topomorphic extension of its MEF.
\end{conj}

Further, a classical result due to Auslander and Yorke in
\cite{AuslanderYorke1980IntervalMaps} states that a minimal tds is either
equicontinuous or sensitive (has sensitive dependence on initial conditions). As
shown by Li, Tu and Ye in \cite{LiTuYe2015MeanSensitivity}, an analogue holds for
the mean versions of these notions as well. Following \cite{LiYeYu2022EquicontinuityAndSensitivityInMeanForms}, 
we say $(X,\varphi)$ is {\em mean $m$-sensitive} if there exists $\eps>0$ such that for
every open set $U\ssq X$ there exist points $x_1\ld x_m\in B_\delta(x)$ with
$\Dmbar(x_1, x_2 \ld x_m)\geq \eps$.
\begin{thm} \label{t.auslander-yorke}
  A minimal tds $(X,\varphi)$ is mean $m$-equicontinuous if and only if it is
  not mean $m$-sensitive.
\end{thm}

It is noteworthy that the finite topomorphic extension structure defined above
is shared by a broad scope of classical examples. Probably the best-know case of
2:1 topomorphic extensions (over the dyadic odometer) are subshifts induced by
the Thue-Morse substitution and their generalisations
\cite{Keane1968GeneralisedMorse}. As Keane already remarked in
\cite{Keane1968GeneralisedMorse}, similar substitutions on alphabets with $m$
symbols should lead to $m$:1 topomorphic extensions in an analogous
way. Moreover, all constant length substitution induce subshifts that have a
finite topomorphic extension structure
\cite{Kamae1972TopologicalInvariantForSubstitutionSystems,Dekking1977ConstantLengthSubstitutionSpectra}. Further
examples include certain irregular Toeplitz flows constructed by Williams
\cite{Williams1984ToeplitzFlows} and by Iwanik and Lacroix
\cite{IwanikLacroix1994NonRegularToeplitz}, and similar examples have been be
obtained in the class of irregular model sets
\cite{FuhrmannGlasnerJaegerOertel2021TameImpliesRegular}. Smooth examples of
minimal skew-products on the two-torus with the same extension structure have
recently been constructed in \cite{HauptJaeger2023IsomorphicExtensions}.\medskip

The paper is organised as follows. In Section \ref{Preliminaries}, we provide
the necessary background and preliminaries. In Section~\ref{FTE}, we discuss
basic properties and alternative definitions of finite-to-one topomorphic
extensions.  The Auslander-Yorke type dichotomy
(Theorem~\ref{t.auslander-yorke}) is shown in Section~\ref{Auslander-Yorke},
alongside with pointwise characterisations of multivariate mean equicontinuity
and sensitivity. The proof of
Theorem~\ref{t.multivariate_mean_equicont_for_topomorphic_extensions} is then
given in Section~\ref{MainResult}, where we also include a brief discussion of
Conjecture~\ref{conj.topomorphic_extensions}.


\section{Notation and preliminaries}\label{Preliminaries}

We assume that the reader is familiar with standard notions of topological
dynamics and ergodic theory, as provided, for instance, in
\cite{Walters1982ErgodicTheory,Petersen1983ErgodicTheory,katok/hasselblatt:1997,Einsiedler2010-ln}
A {\em measure-preserving dynamical system (mpds)} is a quadruple
$(X,\cB,\mu,\varphi)$ consisting of a probability space $(X,\cA,\mu)$ and a
bi-measurable bijective transformation $\varphi:X\to X$ which preserves the
measure $\mu$. In our context, $X$ will in most cases be a compact metric space,
$\cB=\cB(X)$ the Borel $\sigma$-algebra generated by the topology on $X$ and
$\varphi$ a homeomorphism of $X$. Since any compact metric space is Polish,
$(X,\cB(X))$ is a standard Borel space in this situation.\foot{Recall that a
  {\em standard Borel space} is a measurable space of the form $(X,\cB(X))$,
  where $X$ is a Polish space.}

Any measure that is mentioned in the following will implicitly be understood to
be a probability measure, unless explicitly stated otherwise. If $\mu$ is a
measure on a Borel space $(X,\cB(X))$, we denote by
\[
\supp(\mu)\ = \ \{x\in X\mid \mu(U)>0 \textrm{ for all open neighbourhoods } U
\textrm{ of } x\}
\]
the topological support of $\mu$.\smallskip


Suppose now that $(X,\cB,\mu,\varphi)$ and $(Y,\cA,\nu,\psi)$ are two mpds. Then
a measurable map $\pi:X\to Y$ is called a {\em (measure-theoretic) factor map}
and $(Y,\cA,\nu,\psi)$ a {\em (measure-theoretic) factor} of
$(X,\cB,\mu,\varphi)$ if $\pi\circ \varphi=\psi\circ\pi$ holds $\mu$-almost
surely. If additionally there exist subsets $X_0\ssq X,\ Y_0\ssq Y$ of full measure such that
$\pi:X_0\to Y_0$ is a bi-measurable bijection, we call $\pi$ an {\em
  isomorphism} of mpds and say the two systems are {\em (measure-theoretically)
  isomorphic}.

\begin{thm}[Rokhlin's skew product theorem,
    {\cite[Theorem 3.18]{Glasner2003Joinings}}] \label{t.rokhlin} Suppose that
  $(X,\cB)$ and $(Y,\cA)$ are standard Borel spaces, $(X,\cB,\mu,\varphi)$ is an
  ergodic mpds and $(Y,\cA,\mu,\psi)$ is a factor mpds with factor map $\pi:X\to
  Y$. Then there exists a standard Lebesgue space\foot{A {\em standard Lebesgue
      space} $(Z,\cC,\lambda)$ is a standard Borel space $(Z,\cC)$ equipped with
    a probability measure $\lambda$ such that $\cC$ is complete with respect to
    $\lambda$ (that is, any subset of a null set in $\cB$ is again contained in
    $\cC$).} $(Z,\cC,\lambda)$ and an $\cA\otimes\cC$-bi-measurable bijection
  $\rho:Y\times Z\to Y\times Z$ preserving $\nu\otimes\lambda$ such that the two
  systems $(X,\cB,\mu,\varphi)$ and $(Y\times
  Z,\cA\otimes\cC,\nu\otimes\lambda,\rho)$ are isomorphic.

  Further, the transformation $\rho$ can be chosen such that it has skew
  product form
  \[
  \rho : Y\times Z\to Y\times Z \quad , \quad (y,z)\mapsto (\psi(y),\rho_y(z))
  \]
  where $\rho_y:Z\to Z$ preserves the measure $\lambda$ for all $y\in
  Y$. Moreover, the isomorphism $\iota:X\to Y\times Z$ can be chosen such that
  it satisfies $p_Y\circ \iota=\pi$, where $p_Y:Y\times Z\to Y$ denotes the
  canonical projection to $Y$.
\end{thm}

Whenever we invoke Rokhlin's Theorem in the following, we always assume that the
transformation $\rho$ and the isomorphism $\iota$ satisfy the additional
assertions stated above.\smallskip

We will need to use a direct implication of this statement for the structure of
ergodic measures of (measure-theoretically) finite-to-one extensions. Recall
that if $X,Y$ are Polish spaces, $\mu$ is a Borel probability measure on $X$ and
$\nu=\pi_*\mu$, then there exists a mapping
\[
 Y\times \cB(X) \ \to \ [0,1] \quad , \quad (y,A)\mapsto \mu_y(A)
 \]
 such that
 \begin{itemize}
 \item for every $y\in Y$, the mapping $A\mapsto \mu_y(A)$ is a Borel
   probability measure on $X$;
 \item for every $A\in \cB(X)$, the function $y\to \mu_y(A)$ is integrable and
   \[
   \mu(A) \ = \ \int_Y \mu_y(A)\ d\nu(y) \ .
   \]
 \end{itemize}
Such a mapping, which we will also denote as $(\mu_y)_{y\in Y}$, is called a
{\em conditional probability distribution} of $\mu$ over $\pi$. It is unique
modulo modifications on a subset of $Y$ of measure zero and we have
$\mu_y(\pi^{-1}{y})=1$ for $\nu$-almost every $y\in Y$. We refer to the measures
$\mu_y$ as {\em fibre measures}.

Now, let $(X,\cB,\mu,\varphi)$ be an extension of $(Y,\cA,\nu,\psi)$ with
corresponding factor map $\pi$ and suppose that $\mu$ is ergodic. Further,
assume that the support $\supp(\mu_y)$ has finite cardinality $\nu$-almost
surely. Then, since $\mu_y$ is the push-forward of $\lambda$ under the mapping
$z\mapsto \iota^{-1}(y,z)$, the measure $\lambda$ in Theorem~\ref{t.rokhlin} is
supported on a finite set as well. We can therefore assume that the space $Z$
itself is finite, that is, $Z=\{z_1\ld z_k\}$ for some $k\in\N$.  If we let
$\gamma_i(y)=\iota^{-1}(y,z_i)$, $i=1\ld k$, then the situation we arrive at is
the following: We have a measurable multivalued function
\[
\gamma:Y\to X^k\quad ,\quad y\mapsto \gamma(y)= (\gamma_1(y)\ld \gamma_k(y))
\]
 such that
\begin{itemize}
\item $\pi(\gamma_i(y))=y$ holds for $\nu$-almost every $y\in Y$ and all $i=1\ld
  k$;
\item $\gamma_i(y)\neq \gamma_j(y)$ holds for $\nu$-almost every $y\in Y$ and
  all $1\leq i<j\leq k$.
\end{itemize}
and the measure $\mu_\gamma$ is of the form
\begin{equation}\label{e.graph_measure}
\mu_\gamma\ =\ \ktel \sum_{i=1}^k\gamma_{i*}\nu \ .
\end{equation}
We will call a measure of the form (\ref{e.graph_measure}) a {\em graph measure
  of multiplicity $k$}. The term is motivated by the situation where
$(X,\varphi)$ is a skew product system, that is, $X=Y\times \Xi$ and
$\varphi(y,\xi)=(\psi(y),\varphi_y(\xi))$ for all $(y,\xi)\in X$, and $\pi$ is
just the canonical projection from $X=Y\times\Xi$ to $Y$.  In this case
$\gamma_i(y)=(y,\hat\gamma_i(y))$ for some measurable function
$\hat\gamma_i:Y\to\Xi$ and the measure $\mu_\gamma$ is supported on the union of
the graphs of $\hat\gamma_1\ld \hat\gamma_k$.\smallskip

Altogether, we obtain the following.
\begin{lem}\label{l.graph_measures}
  Suppose that the mpds $(X,\cB,\mu,\varphi)$ is an extension of
  $(Y,\cA,\nu,\psi)$. Further, assume that $\mu$ is ergodic and
  $\sharp\supp(\mu_y)<\infty$ $\nu$-almost surely.  Then $\mu$ is a graph
  measure of finite multiplicity $k\in\N$. In particular, we have
  $\sharp\supp(\mu_y)=k$ for $\nu$-almost every $y\in Y$.
\end{lem}

Finally, given a tds $(X,\varphi)$, a $\varphi$-invariant measure $\mu$ and
$n\in \N$, an $n$-fold (self-)joining of $\mu$ is a $\varphi^{\times
  n}$-invariant measure $\hat\mu$ on $X^n$ that satisfies $\pi_{i*}\hat\mu=\mu$
for $i=0\ld n$. Here
\[
\varphi^{\times n} :X^n\to X^n \quad , \quad (x_1\ld x_n)\mapsto
(\varphi(x_1)\ld \varphi(x_n)) \ .
\]
If $(Y,\psi)$ is factor of $(X,\varphi)$ with factor map $\pi:X\to Y$ and $\hat
\mu$ is supported on the set
\[
X^n_\pi \ = \ \{x\in X^n\mid \pi(x_1)=\ldots = \pi(x_n)\} \ ,
\]
then $\hat\mu$ is called an {\em $n$-fold joining over the common factor
  $(Y,\psi)$}. There is always at least one such joining: if $\mu$ disintegrates
as $(\mu_y)_{y\in Y}$, then it can be obtained by integration of the fibre
measures $\hat\mu_y = \bigotimes_{i=1}^n \mu_y$ with respect to $\nu$. We refer
to \cite{Glasner2003Joinings,DeLaRue2020Joinings} for further background on
joinings.


\section{Finite topomorphic extensions} \label{FTE}

The aim of this section is to provide two seemingly weaker, but alternative
definitions of finite topomorphic extensions. First, instead of requiring the
existence of an $m$-valued mapping $\gamma:Y\to X^m$, $m\in\N$, whose
corresponding point set supports all $\varphi$-invariant measures, it suffices to
require that all these measures are supported on some measurable set that
intersects the fibres $\pi^{-1}(y)$ in at most $m$ points. The {\em `graph
  structure'} of this set then comes for free. Secondly, it also suffices to
require that any $\varphi$-invariant measure is supported on at most $m$ points
in every fibre. In both cases, the proof of the equivalence to the original
definition hinges on an application of Rokhlin's Skew Product
Theorem. \smallskip

\begin{prop} \label{p.topomorphic_extensions_equivalent_characterisations}
  Suppose $(X,\varphi)$ is a tds with uniquely ergodic MEF $(Y,\psi)$ and
  corresponding factor map $\pi:X\to Y$. Denote by $\nu$ the unique
  $\psi$-invariant measure.  Then the following are equivalent for all $m\in\N$.
  \romanlist
\item $(X,\varphi)$ is an $m$:1 topomorphic extension of $(Y,\psi)$.
\item $m$ is the least integer such that there exists a measurable set $M\ssq X$
  with $\mu(M)=1$ for all $\varphi$-invariant Borel probability measures on $X$
  and $\sharp\pi^{-1}\{y\}\cap M=m$ for $\nu$-almost every $y\in Y$.
\item For every $\varphi$-invariant measure $\mu$ on $X$, the fibre measures
  $\mu_y$ are $\nu$-almost surely supported on $m$ points.
  \listend
\end{prop}
In fact, the assumption of unique ergodicity of $(Y,\psi)$ is not strictly
necessary. However, since this is always satisfied for minimal systems and these
are our main focus, we use it here both for the sake of simplicity and because
some additional information is available in this case (see
Addendum~\ref{p.topomorphic_extensions_measures}).

\proof\quad Let the least integers for which the statements in (i), (ii) and
(iii) hold be denoted by $m$, $m'$ and $m''$, respectively. A priori, we also
allow values $+\infty$. However, we show $m\geq m'\geq m'' \geq m$, which proves
the asserted equivalence.
\smallskip

    {$m\geq m'$}: suppose $\gamma:Y\to X^m$ is such that the corresponding point
    set $\Gamma$ supports all $\varphi$-invariant probability measures. Then it
    suffices to set $M=\Gamma$.
    \smallskip
    
	{$m'\geq m''$}: Let $M_y = M\cap \pi^{-1}\{y\}$.
	Then $\mu(M)=1$ implies $\mu_y(M_y)=1$ for $\nu$-almost
	every $y\in Y$. Hence, $\mu_y$ is supported on $m'$ points $\nu$-almost surely.
	\smallskip

	{$m''\geq m$}: We first show that there exist at most $m$ different ergodic
	$\varphi$-invariant measures on $X$. Assume otherwise and let $\mu^1\ld
	\mu^{m+1}$ denote different $\varphi$-invariant ergodic measures. For each of
	these measures, the fibre measures $\mu^i_y$ are supported on at most $m''$
	points $\nu$-almost surely. This means that if we apply Rokhlin's Skew Product
	Theorem to $\mu^i$, we obtain that $(X,\cB(X),\mu^i,\varphi)$ is isomorphic via
	an isomorphism $\iota_i$ to a skew product system $(Y\times Z_i,\cB(Y)\otimes
	\cC_i,\nu\otimes\lambda_i,\rho_i)$ where $\lambda_i$ is the equidistribution (by
	ergodicity) on a finite set $\hat Z_i=\{z^i_1\ld z^i_{k_i}\}$ with $k_i\leq
	m''$. The measure $\mu^i_y$ is then supported on the set $\hat
	X_i(y)=\iota_i^{-1}\kl\{y\}\times\hat Z_i\kr$. By ergodicity, the sets
	$X_1(y)\ld X_{m''+1}(y)$ are $\nu$-almost surely disjoint. However, this means
	that the fibre measures $\mu_y$ of the measure
	$\mu=\frac{1}{m''+1}\sum_{i=1}^{m''+1} \mu^i$ are supported on the sets
	$X(y)=\bigcup_{i=1}^{m''+1}X_i(y)$ of cardinality $\sum_{i=1}^{m''+1}k_i > m''$,
	contradicting the assumption.
	\smallskip

Hence, if $\mu^1\ld \mu^\ell$ are all the different ergodic $\varphi$-invariant
measures on $X$, then $\ell\leq m''$. More precisely, the above argument yields
that if the fibre measures of $\mu^i$ are supported on $k_i$ points in
$\nu$-almost every fibre, then $\tilde m=\sum_{i=1}^\ell k_i \leq m''$. Given
$j\in\{1\ld \tilde m\}$, we can write $j$ in a unique way as $j=\sum_{i=1}^{s-1}
k_i+t$, where $s\in\{1\ld \ell\}$ and $t\in\{1\ld k_s\}$. If we then define
\[
\gamma_j(y) \ = \ \iota^{-1}_s(z^s_t)
\]
we obtain a measurable mapping $\gamma:Y\to X^{\tilde m}$ that satisfies
$\mu^i(\Gamma)=1$ for all $i=1\ld \ell$, and hence $\mu(\Gamma)=1$ for all
$\varphi$-invariant measures by ergodic decomposition.
\qed\medskip

In combination with Lemma~\ref{l.graph_measures}, the preceding proof of
Proposition~\ref{p.topomorphic_extensions_equivalent_characterisations} actually
yields some additional information.

\begin{add}\label{p.topomorphic_extensions_measures}
  Suppose we are in the situation of
  Proposition~\ref{p.topomorphic_extensions_equivalent_characterisations}. Then
  there exist at most $m$ $\varphi$-invariant measures $\mu_1\ld \mu_\ell$,
  where $\ell\leq m$, all of which are graph measures. Further, if $k_i$ denotes
  the multiplicity of the graph measure $\mu_i$, then $\sum_{i=1}^\ell k_i=m$.
\end{add}


\section{An Auslander-Yorke type dichotomy} \label{Auslander-Yorke}

Before we turn to a closer inspection of multivariate mean equicontinuity, we
first want to collect some basic facts concerning the mappings
$\Dmmax$, $\Dm$ and $\Dmbar$ introduced
above. We will refer to these as {\em multidistances}.

\subsection{The multidistances $\Dm$ and $\Dmmax$}

Given $m\in\N$ and $\mathbf{x}\in X^m$, we define $\Rep{\mathbf{x}}{i}{z}$ by $
\Rep{\mathbf{x}}{i}{z}_j = x_j$ if $j\neq i$ and $\Rep{\mathbf{x}}{i}{z}_i=z$,
i.e. $\Rep{\mathbf{x}}{i}{z}=(x_1,\ldots,x_{i-1},z,x_{i+1}\ldots x_m)$.
We say that a map $\mathrm{D}: X^m \to \R^+_0$ satisfies the \textdef{polygon
  inequality} if and only if for any $\mathbf{x}\in X^m$ and $z \in X$ we have
\[
\mathrm{D}(\mathbf{x}) \leqslant \sum_{i=1}^m \mathrm{D}\Rep{\mathbf{x}}{i}{z}
\ .
\]
If additionally $\mathrm{D}$ is \textdef{symmetric}, that is, for any tuple
$(x_1,\ldots,x_m) \in X^{m+1}$ and any permutation $\sigma \in \Sym(m)$ we have
\[
\mathrm{D}(x_1,\ldots,x_m) = \mathrm{D}(x_{\sigma(1)},\ldots,x_{\sigma(m)})  \ ,
\]
and \textdef{positive semi-definite}, meaning $\mathrm{D}(x,\ldots,x) = 0$ for
all $x\in X$, we call $\mathrm{D}$ an \textdef{$m$-multidistance}.

\begin{rem}
	A notable example of a multidistance, not needed in this work, is the
        measure of the simplex spanned by $m$ points in euclidean space
        $\R^{m-1}$.
\end{rem}

\begin{lem}
	$\Dm$ and $\Dmmax$ are $m$-multidistances.
\end{lem}
\proof\quad
	Symmetry and positive semi-definiteness are clear.  We will show the
        polygon inequality.  For $m=1$ we have the original metric and the
        original triangle inequality so let $m \geqslant 2$.  Let
        $\mathbf{x}=(x_1,\ldots,x_m) \in X^{m}$ and $z\in X$ be arbitrary.
		\smallskip

	First we consider $\Dmmax$.  Let $a,b \in \ml1,\ldots,m\mr$ be
        maximizing indices, that is, $ d(x_a,x_b)=D^{\max{}}_m(x_1,\ldots,x_m)
        $.  Then we have
	\[ D^{\max{}}_m \Rep{\mathbf{x}}{i}{z} \geqslant D^{\max{}}_m(\mathbf{x}) \]
	for all $i \notin \ml a,b\mr$.  So a fortiori we have
	\[ \summe i0m D^{\max{}}_m \Rep{\mathbf{x}}{i}{z}
        \geqslant D^{\max{}}_m(\mathbf{x}) \,. \]

	Now we consider $\Dm$.  Pick minimizing indices $a,b \in \ml
        1,\ldots,m\mr$, that is, $d(x_a,x_b) = D_m(x_1,\ldots,x_m)$.  For ease
        of notation, let $x_{m+1} = z$.  Likewise, for any $i \in \ml
        1,\ldots,m\mr$, pick indices $a_i,b_i \in \ml 1,\ldots,m+1\mr\setminus
        \ml i \mr$ such that $d(x_{a_i},x_{b_i}) =
        D_m(\Rep{\mathbf{x}}{i}{x_{m+1}})$.  If for some $j \in \ml 1,\ldots,m
        \mr$ neither $a_j=m+1$ nor $b_j=m+1$, we are done, as then
	\[
        \sum_{i=1}^m D_m \Rep{\mathbf{x}}{i}{x_{m+1}} \geqslant
        D_m(\Rep{\mathbf{x}}{j}{x_{m+1}}) = d(x_{a_j},x_{b_j}) \geqslant
        D_m(\mathbf{x}) \ .
        \]
        So without loss of generality $b_i = m+1$ for any $i$.  Now observe that
        $i \mapsto a_i \in \ml 1,\ldots, m\mr$ can not be constant, as $a_i \neq
        i$ by construction.  Pick $k,l\in \ml 1, \ldots, m \mr$ such that $a_k
        \neq a_l$.  Now the triangle inequality implies
	\begin{align*}
		\sum_{i=1}^m D_m\Rep{\mathbf{x}}{i}{x_{m+1}} &\geqslant D_m
                \Rep{\mathbf{x}}{k}{x_{m+1}} + D_m \Rep{\mathbf{x}}{l}{x_{m+1}}
                \\ &= d(x_{a_k},x) + d(x_{a_l},x) \\ &\geqslant
                d(x_{a_k},x_{a_l}) \geqslant D_m(\mathbf{x}) \,. \qedhere
	\end{align*}
        \medskip

        The multidistance property of $\Dm$ directly carries over to $\Dmbar$ in
        the limit.
        \begin{cor}
          If $(X,\varphi)$ is a tds, then $\Dmbar$ is
          an $m$-multidistance
        \end{cor}

        \subsection{Pointwise multivariate mean equicontinuity and sensitivity}

        Recall that a tds $(X,\varphi)$ is mean $m$-equicontinuous, with $m\geq
        2$, if for every $\eps>0$ there exists $\delta>0$ such that
        $\Dmmax(\mathbf {x})<\delta$ implies $\Dmbar(\mathbf{x})<\eps$ for all
        $\mathbf{x}\in X^m$. Conversely, $(X,\varphi)$ is mean $m$-sensitive if
        for some $\eps>0$ and every open set $U\ssq X$ there exists
        $\mathbf{x}\in U^m$ with $\Dmbar(\mathbf{x})\geq \eps$.  The aim of this
        section is to provide pointwise characterisations of these
        properties. Note that neither minimality nor any other recurrence
        assumption is required here.\smallskip

         Given $\eps>0$, a point $x \in X$ is called an \textdef{mean
           $\eps$-$m$-equicontinuity point} if and only if there is a 
		   $\delta>0$ such that $\Dmbar(x,x_2,\ldots,x_m)  < \varepsilon$ for all $x_2,\ldots, x_m\in B_\delta(x)$. It is called a
         \textdef{mean $m$-equicontinuity point} if it is a mean
         $\eps$-$m$-equicontinuity point for all $\eps>0$. Conversely, we say
         point $x\in X$ is {\em mean $m$-sensitive (mean $\eps$-$m$-sensitive)}
         if it is not a mean $m$-equicontinuity point (mean
         $\eps$-$m$-equicontinuity point).

         \begin{prop}\label{p.pointwise_mean_equicontinuity}
           Let $(X,\varphi)$ be a tds and $m\geq 2$.  \alphlist
         \item $(X,\varphi)$ is mean $m$-equicontinuous if and only if every
           $x\in X$ is a mean $m$-equicontinuity point. 
         \item $(X,\varphi)$ is mean $m$-sensitive if and only if there exists
           $\eps>0$ such that every $x\in X$ is mean $\eps$-$m$-sensitive.
           \listend
         \end{prop}
          Note that for the case $m=2$ this result is contained in
          \cite{LiTuYe2015MeanSensitivity}. In the following proof, we
          essentially adapt the respective arguments to the multivariate case.
          \proof \alphlist
         \item If $(X,\varphi)$ is mean $m$-equicontinuous, then it is obvious
           that all points are mean $m$-equicontinuity points. Conversely,
           assume that every $x\in X$ is a mean $m$-equicontinuity point and fix
           $\eps>0$. For every $x\in X$, choose $\delta(x)>0$ such that
           $\Dmbar(x,x_2\ld x_m)<\frac{\eps}{m}$ for all $x_2\ld x_m\in
           B_{\delta(x)}(x)$. Then $\cU=\{B_{\delta(x)}(x)\mid x\in X\}$ is an
           open cover of $X$. Choose any $\delta>0$ smaller than the Lebesgue
           covering number of $\cU$.  Then, if $\mathbf{x} =(x_1\ld x_m)\in X^m$
           satisfies $\Dmmax(x_1\ld x_m)<\delta$, we have that $\{x_1\ld
           x_m\}\ssq B_\delta(x_1)\ssq B_{\delta(\xi)}(\xi)$ for some $\xi\in
           X$. Consequently, the polygonal inequality yields
           \[
           \Dmbar(\mathbf{x}) \ \leq \ \sum_{i=1}^m
           \Dmbar(\Rep{\mathbf{x}}{i}{\xi}) \ < \ m\cdot \frac{\eps}{m} \ =
           \ \eps \ .
           \]
           This shows the mean $m$-equicontinuity of $(X,\varphi)$.
         \item Suppose that $(X,\varphi)$ is mean $m$-sensitive. Then there
           exists some $\eta > 0$ such that for any non-empty open subset $U
           \subseteq X$ there is $\mathbf{x} \in U^m$ with $\Dmbar(\mathbf{x}) >
           \eta$.  Let $x \in X$, $\delta>0$ and $U = B_\delta(x)$. Choose
           $\mathbf{x} \in U^{m}$ with $\Dmbar(\mathbf{x}) > \eta$.  Then by the
           multi-triangle inequality we have
\[ \sum_{i=0} \Dmbar\kl \Rep{\mathbf{x}}{i}{x}
\kr \geqslant \Dmbar(\mathbf{x}) > \eta \,.
\]
So $\Dmbar\kl \Rep{\mathbf{x}}{i}{x} \kr > \frac{\eta}{m}$ for at least for one
$i \in \ml 0,\ldots, m \mr$. This yields that $x$ is a mean $\eps$-$m$-sensitive
point with $\eps=\frac{\eta}{m}$. As $x\in X$ was arbitrary, this proves the
first implication. The converse direction is again obvious.\qed\medskip
\listend

\subsection{Auslander-Yorke type dichotomy}

We again refer to \cite{LiTuYe2015MeanSensitivity} for analogous results on the
case $m=2$. The key observation is the following.

\begin{lem}
	\label{cor:almostauslanderyorkedichotomy}
	Let $(X,\varphi)$ be transitive with transitivity point $x$.
	Then either $x$ is a mean $m$-equicontinuity point
	or the system is mean $(m+1)$-sensitive.
\end{lem}
\begin{proof}
        Suppose that $x$ is a transitive point, but not a mean
        $m$-equicontinuity point. Let $U\ssq X$ be open. By definition, $x$ is
        mean $\varepsilon$-$m$-sensitive for some $\varepsilon > 0$.  However,
        this means that $\varphi^n(x)$ is mean $\varepsilon$-$m$-sensitive for
        any $n \in \N$.  As $x$ is a transitivity point, $\varphi^n(x) \in U$
        for some $n\in\N$. Since $U$ is open, we can choose $x_2\ld x_m\in U$
        such that $\Dmbar(\varphi^n(x),x_2\ld x_m)\geq \eps$. This proves mean
        $m$-sensitivity of $(X,\varphi)$.
\end{proof}

For a minimal system $(X,\varphi)$, this means that either the system is mean
$m$-sensitive or all points are mean $m$-equicontinuity points. Due to
Proposition~\ref{p.pointwise_mean_equicontinuity}(a), the latter implies that
$(X,\varphi)$ is mean $m$-equicontinuous. This proves
Theorem~\ref{t.auslander-yorke}, which we restate here as
\begin{cor}[Auslander-Yorke type Dichotomy]
	A minimal tds $(X,\varphi)$ is mean $m$-equi\-continuous if and only if
        it is not mean $m$-sensitive.
\end{cor}

\subsection{Non-continuity of the Besicovitch multidistance}

We want to close this section by pointing out an important difference between
mean equicontinuity and multivariate mean equicontinuity ($m>2$). As mentioned
in the introduction, mean equicontinuity can be defined as the continuity of the
Besicovitch pseudo-metric. This is different in the multivariate case, since
continuity of the Besicovitch $m$-distance equally implies mean
equicontinuity. Hence, if $(X,\varphi)$ is mean $m$-equicontinuous for some
$m>2$, but not mean equicontinuous, the corresponding Besicovitch $m$-distance
cannot be continuous. Examples of this type are given, for instance, by
Thue-Morse subshifts (\cite{Keane1968GeneralisedMorse}, compare \cite[Theorem
  4.6]{LiYeYu2022EquicontinuityAndSensitivityInMeanForms} and its proof) or by
certain irregular Toeplitz flows
\cite{Williams1984ToeplitzFlows,IwanikLacroix1994NonRegularToeplitz}.

\begin{lem}
	\label{lem:discontinuityofDmbar}
	Assume that $(X,\varphi)$ has an infinite MEF.
	If $\Dmbar$ is continuous, then $(X,\varphi)$ is mean
        equicontinuous.
\end{lem}
\begin{proof}
	Let $\pi: X \to Y$ be the factor map to the MEF
	$(Y,\psi)$.  Pick a $\psi$-invariant metric $d_Y$ on $Y$ and note that
	we can switch to equivalent metrics in order to ensure that $d(x,y) >
	d_Y(\pi(x),\pi(y))$ for any $x,y \in X$.  Pick any $x_2 \in X$ and
	choose $x_3,\ldots,x_{m} \in X$ such that $\pi(x_i) \neq \pi(x_j)$ for any
	$2\leq i < j \leq m$.  Define $c=\min_{2 \leqslant i < j \leqslant m}
	d_Y(\pi(x_i),\pi(x_j)) > 0$.  Now note that if $d_Y(\pi(x),\pi(x_2)) <
	\frac{c}{2}$, then $\min_{3 \leqslant i \leqslant m}
	d_Y(\pi(x),\pi(x_j)) > \frac{c}{2}$.  By invariance of $d_Y$ this holds
	for any iterate.

	Observe that $\Dmbar(x_2,x_2\ldots,x_m) = 0$, so that continuity of
        $\Dmbar$ implies
	\begin{align*}
	 	\limsup_{n\to\infty} \frac 1n \summe k1n \Dm(\varphi^k(x),\varphi^k(x_2),\ldots, \varphi^k(x_{m})) =
                \Dmbar(x,x_2\ldots,x_{m}) \ \xrightarrow{x\to x_2} \ 0
                \ . 
	\end{align*}
	Expanding we obtain
	\begin{eqnarray}
		\lefteqn{\nonumber \Dm(\varphi^k(x),\varphi^k(x_2),\ldots,
                  \varphi^k(x_{m}))} \notag \\ \label{eq:expandedDm1} & =
                &\min\Big\{ d\kl \varphi^k(x), \varphi^k(x_2) \kr, \min_{3
                  \leqslant j \leqslant m} d\kl \varphi^k(x),\varphi^k(x_j) \kr,
                \min_{2 \leqslant i < j \leqslant m} d\kl
                \varphi^k(x_i),\varphi^k(x_j)) \kr \Big\} \,.
	\end{eqnarray}
	Note that the terms in \refe{expandedDm1} are uniformly bounded from
        below by $\frac c2$.  So it is the term $d\kl
        \varphi^k(x),\varphi^k(x_j) \kr$ that must be responsible for
        $\Dmbar(x,x_2,\ldots,x_m)$ getting arbitrarily small for $x$
        sufficiently close to $x_2$.  This in turn implies that
        $d_\textsf{B}(x,x_2)$ gets arbitrarily small for $x$ sufficiently close
        to $x_2$.  So $x_2$ is a mean equicontinuity point.

	As $x_1 \in X$ was arbitrary, this shows that all points in $X$ are mean
        equicontinuity points. Hence, by
        Proposition~\ref{p.pointwise_mean_equicontinuity}(a), $(X,\varphi)$ is
        mean equicontinuous.
\end{proof}


\section{Multivariate mean equicontinuity of finite topomorphic extensions}
\label{MainResult}

We again consider a tds $(X,\varphi)$ with MEF $(Y,\psi)$ and corresponding
factor map $\pi:X\to Y$. Given $m\in\N$ and $\mathbf{x}\in X^m$, we let
$\Dmymax(\mathbf{x})=\max_{1\leq i<j\leq m} d_Y(\pi(x_i),\pi(x_j))$, where $d_Y$
denotes the metric on $Y$. Then, we call $(X,\varphi)$ {\em factor mean
  $m$-equicontinuous} if for all $\eps>0$ there exists $\delta>0$ such that
\[
\Dmymax(\mathbf{x})<\delta \quad \textrm{implies} \quad \Dmbar(\mathbf{x})<\eps \ .
\]
Due to the continuity of the factor map $\pi$, it is obvious that factor mean
$m$-equicontinuity implies mean $m$-equicontinuity. Therefore,
Theorem~\ref{t.multivariate_mean_equicont_for_topomorphic_extensions} is a
direct consequence of the following equivalence.
\begin{thm}\label{t.topomorphic_extension_characterisation}
  Let $m\in\N$.  A minimal tds $(X,\varphi)$ is an $m$:1 topomorphic extension
  of its MEF if and only if it is factor mean $(m+1)$-equicontinuous, but not
  factor mean $k$-equicontinuous for any $k\leq m$.
\end{thm}
\proof\quad We first assume that $(X,\varphi)$ is an $m$:1 topomorphic extension
of its MEF, so that all $\varphi$-invariant measures are supported on a set
$\Gamma=\{\gamma_i(y)\mid i=1\ld m,\ y\in Y\}$, where $\gamma:Y\to X^m$ is
measurable. Let the $\nu$ be the unique invariant measure on $(Y,\psi)$.  By
Addendum~\ref{p.topomorphic_extensions_measures}, there exist a finite number
$\ell\leq m$ of ergodic measures on $(X,\varphi)$, which we denote by
$\mu_1,\ldots,\mu_k$.  We assume without loss of generality that the metric
$d_Y$ on $Y$ is $\psi$-invariant and that $d(x_1,x_2) \geqslant
d_Y(\pi(x_1),\pi(x_2))$.
	
Let $\eta,\rho > 0$, where both parameters are arbitrary at the moment, but will
be specified further later.  By Lusin's Theorem, there is a compact set $K =
K_\eta \subseteq Y$ with $\nu(K) > 1-\eta$ such that $\gamma\vert_K$ is
continuous.  Let
	\[
        U \ = \  B_\rho\kl \bigcup_{i=1}^m \gamma_i(K) \kr \ .
        \]
        Then, since all $\varphi$-invariant measures are supported on $\Gamma$
        and project to $\nu$, we have that $\mu(U)>1-\eta$ for any
        $\varphi$-invariant measure $\mu$. As a consequence of the semi-uniform
        ergodic theorem \cite{sturman/stark:2000}, we obtain
	\[
        \liminf_{n \to \infty} \frac 1N \sum_{i=1}^n \ind_U\circ \varphi^n(x)
        \geqslant 1-\eta \] for any $x \in X$. Now, let $\mathbf{x}\in X^{m+1}$
        and consider the set of simultaneous hitting times
	\[
        T(\mathbf{x}) \ = \ \ml t \in \N \mm \varphi^t(x_j) \in U \text{ for all
        } j \in \ml 1, \ldots, m+1 \mr \mr \ .
        \]
	For any $\mathbf{x} = (x_0,\ldots,x_m) \in X^{m+1}$, this set has lower
        asymptotic density
	\[
		a(\mathbf{x}) \ = \ \liminf_{N \to \infty} \frac 1N \#\ml
                T(\mathbf{x})\cap \{ 1, \ldots, N \}\mr \ \geqslant \ 1 - (m+1)
                \eta \ .
	\]

 Now, fix $\varepsilon > 0$.  We need to find $\delta > 0$ such that
 $\Dmiymax(\pi(\mathbf{x})) < \delta$ implies $\Dmibar(\mathbf{x}) < \varepsilon$.
 We have that
	\begin{align*}
		\Dmibar(\mathbf{x}) &\leqslant (1-a(\mathbf{x})) \cdot \diam(X) +
                a(\mathbf{x}) \cdot \kappa(\mathbf{x}) \\ & \leqslant (m+1)
                \cdot \eta \cdot \diam(X) + \kappa(\mathbf{x})
	\end{align*}
	holds for any $\mathbf{x} \in X^{m+1}$, where
	\[
        \kappa(\mathbf{x}) \ = \ \sup_{n\in T(\mathbf{x})} \Dmi(\varphi^n(\mathbf{x}))
        \]
         As $X$ is compact, $\diam(X) < \infty$. We can therefore choose $\eta$
         small enough such that
	\[ (m+1) \cdot \eta \cdot \diam(X) \ < \  \frac{\varepsilon}{2} \ . \]
	It remains to find $\delta > 0$ such that $\Dmiymax(\pi(\mathbf{x})) <
        \delta$ implies $\kappa(\mathbf{x})\leq\frac{\varepsilon}{2}$.  As the
        $\gamma_j$ are continuous on $K$, there is $\alpha > 0$ such that
	\[
        d_Y(y_1,y_2) < \alpha \Longrightarrow
        d(\gamma_j(y_1),\gamma_j(y_2)) < \frac{\varepsilon}{6}
        \]
        holds for all $y_1,y_2 \in K$ and any $j \in \ml 1, \ldots, m\mr$.  We
        now fix $\rho <
        \min\left\{\frac{\alpha}{3},\frac{\varepsilon}{6}\right\}$ and let
        $\delta < \frac{\alpha}{3}$.\smallskip

	Suppose that $\Dmiymax(\pi(\mathbf{x})) < \delta$ and let $n \in
        T(\mathbf{x})$ be arbitrary.  For any $i\in\ml 1,\ldots, m+1\mr$ we have
        $\varphi^n(x_i) \in U$.  So there is $j_i \in \ml 1,\ldots,m \mr$ and
        $y_i \in K$ such that $d\kl \varphi^n(x_i),\gamma_{j_{i}}(y_i) \kr <
        \rho$.  Clearly $d_Y(y_i,\psi^n(\pi(x_i))) < \rho < \frac{\alpha}{3}$.
        Further $d_Y(\psi^n(\pi(x_i)),\psi^n(\pi(x_j))) < \delta <
        \frac{\alpha}{3}$.  Invariance and triangle inequality together imply
        that $d_Y(y_i,y_j) < \alpha$.  So for any $k \in \ml 1,\ldots, m\mr$ we
        have $d(\gamma_k(y_i),\gamma_k(y_j)) < \frac{\varepsilon}{6}$.
	
	Note that we have a map $i \mapsto j_i$ that goes from $\ml 1,\ldots, m+1
        \mr$ to $\ml 1, \ldots, m\mr$.  Therefore, the pigeon hole principle
        implies the existence of indices $a \neq b$ from $\ml 1, \ldots, m+1 \mr$
        such that $j_{a} = j_{b}$.  We write $k = j_{a}$ for that common value
        and obtain
	\begin{align*}
		d(\varphi^n(x_a),\varphi^n(x_b)) &\leqslant d\kl
                \varphi^n(x_a),\gamma_{k}(y_a) \kr + d\kl
                \gamma_{k}(y_a),\gamma_{k}(y_b) \kr + d\kl
                \gamma_{k}(y_b),\varphi^n(x_b) \kr\\ &< \rho +
                \frac{\varepsilon}{6} + \rho < \frac{\varepsilon}{6} +
                \frac{\varepsilon}{6} +\frac{\varepsilon}{6} =
                \frac{\varepsilon}{2} \,.
	\end{align*}
	So for any $n \in T(\mathbf{x})$ there are indices $a \neq b$ such that
        $d(\varphi^n(x_a),\varphi^n(x_b))< \frac{\varepsilon}{2}$ and therefore
	\[ \Dmi(\varphi^n(\mathbf{x})) < \frac{\varepsilon}{2} \,.\]
	This now implies $\kappa(\mathbf{x}) \leq \frac{\varepsilon}{2}$ as
        required and proves the mean $(m+1)$-equicontinuity of
        $(X,\varphi)$.\smallskip

For the converse direction, assume that $(X,\varphi)$ is factor mean
$(m+1)$-equicontinuous and suppose for a contradiction that it is not a $k$:1
topomorphic extension of its MEF $(Y,\psi)$ for some $k\leq m$. Then, by
Proposition~\ref{p.topomorphic_extensions_equivalent_characterisations}(iii),
there exists some $\varphi$-invariant measure $\mu$ such that
$\sharp\supp(\mu)\geq m+1$ $\nu$-almost surely.  If we consider the $(m+1)$-fold
joining $\hat \mu$ over the factor $\pi$ given by
$\hat\mu_y=\bigotimes_{i=1}^{m+1}\mu_y$, we have that
\[
        \hat\mu\kl X^{m+1}\smin W_{m+1}\kr\ > \ 0 \ ,
\]
where $W_{m+1}=\{\mathbf{x}\in X^{m+1}\mid \exists i\neq j: x_i=x_j\}$.

Due to the Ergodic Decomposition Theorem, there also exist an ergodic joining
$\tilde \mu$ with this property. However, as $\Dmi(\mathbf{x})>0$ on
$X^{m+1}\smin W_{m+1}$, we obtain $\int_{X^{m+1}} \Dmi(\mathbf{x})
\ d\tilde\mu(\mathbf{x}) >0$ and hence $\Dmibar(\mathbf{x})>0$
$\tilde\mu$-almost surely. This shows the existence of points
$(x_1,\ldots,x_{m+1}) \in X^{m+1}$ with the property that $\pi(x_1) = \ldots =
\pi(x_{m+1})$ and $\Dmibar(x_1,\ldots,x_{m+1}) > 0$, contradicting the factor
mean $(m+1)$-equicontinuity.  \qed\medskip
        
\begin{rem}
  The above proof follows the same overall strategy as the proof of the
  equivalence of mean equicontinuity and a topo-isomorphic extension structure
  in \cite{DownarowiczGlasner2015IsomorphicExtensionsAndMeanEquicontinuity},
  with the necessary modifications for the multivariate case. However, the
  statement of Theorem~\ref{t.topomorphic_extension_characterisation} is weaker
  in the sense that we have to replace mean $m$-equicontinuity by factor mean
  $m$-equicontinuity. The reason for this is the following.

  In \cite{DownarowiczGlasner2015IsomorphicExtensionsAndMeanEquicontinuity}, the
  direction from mean equicontinuity to the extension structure is proved by a
  contradiction argument, showing that any system that is not topo-isomorphic to
  its MEF is mean sensitive. The core part of the argument is to show the
  existence of a tuple $\mathbf{x}\in X^{m+1}$, located in some fibre
  $\pi^{-1}(y)$, such that $\Dmbar(\mathbf{x})>0$. This is done in a similar way
  in the proof above. What would be missing in order to show mean
  $m$-sensitivity (and not just lack of factor mean $(m+1)$-equicontinuity) is to
  prove that this tuple $\mathbf{x}$ can be found within an arbitrarily small
  ball. However, in the case $m=1$ this comes for free, since it is known that
  $\pi(x_1)=\pi(x_2)$ implies $d_\textsf{B}(x,y)=0$ and therefore $\inf_{n\in\Z}
  d(\varphi^n(x_1),\varphi^n(x_2))=0$ in mean equicontinuous systems. Thus, one
  may simply replace $x_1$ and $x_2$ by suitable iterates in order to bring the
  two Besicovitch-separated points arbitrarily close.

  When $m\geq2$, however, the analogous statement is not true anymore. The fact
  that $\pi(x_1)=\ldots=\pi(x_{m+1})$ does not imply $\inf_{n\in\N}
  \Dmmax(\varphi^n(\mathbf{x}))=0$. In fact, a careful analysis of the examples
  in \cite{HauptJaeger2023IsomorphicExtensions} reveals that the $m$:1
  topomorphic extensions of irrational rotations constructed there allow no
  $(m+1)$-tuples in a single fibre that satisfy $\inf_{n\in\N}
  \Dmimax(\varphi^n(\mathbf{x}))=0$ and $\Dmibar(\mathbf{x}))=0$ at the same
  time.

  Altogether, this means that any argument allowing to prove
  Conjecture~\ref{conj.topomorphic_extensions} needs to be substantially
  different and must take into account the behaviour of points across different
  fibres. We leave this problem open here.
\end{rem}

\begin{rem}
  Given Theorem~\ref{MainResult}, the validity of
  Conjecture~\ref{conj.topomorphic_extensions} would imply that in the minimal
  case factor mean $m$-equicontinuity and mean $m$-equicontinuity are
  equivalent. We note that this is not true in the general (non-minimal)
  case. Simple counterexamples are given by Morse-Smale systems on the unit
  interval: suppose that $f:[0,1]\to[0,1]$ is a homeomophism with a finite
  number of fixed points $0=x_0<x_1<\ldots < x_k=1$ and $\nLim f^n(x)=x_i$ holds for
  all $x\in (x_{i-1},x_i]$ and $i=1\ld k$. Then it is easy to check
    that the MEF of $([0,1],f)$ is trivial (a singleton $Y_0=\{0\}$), the system
    is a $k$:1 topomorphic extension of $Y_0$, but it is always mean
    $3$-equicontinuous.
\end{rem}



\end{document}